\pdfoutput=1
\documentclass[colorlinks,11pt,reqno]{amsart}

\usepackage{verbatim}
\usepackage{times}
\usepackage[a4paper]{geometry}
\usepackage{tabls}
\usepackage{url}
\usepackage[linktocpage = true]{hyperref}
\hypersetup{
 colorlinks = true,
 citecolor = blue,
}
\usepackage[draft]{fixme}
\fxsetup{layout=footnote, marginclue}
\usepackage{color}
\definecolor{red}{rgb}{1,0,0}
\definecolor{green}{rgb}{0,1,0}
\definecolor{blue}{rgb}{0,0,1}
\definecolor{refkey}{gray}{.625}
\definecolor{labelkey}{gray}{.625}

\usepackage[lite]{amsrefs}
\usepackage{amsfonts}
\usepackage{amssymb}
\usepackage{mathrsfs}
\usepackage{amsmath}
\usepackage{amsthm}
\usepackage{amscd}
\usepackage{enumerate}
\usepackage{paralist}
\usepackage{xypic}
\usepackage{graphicx}
\usepackage{mathtools}
\usepackage[all,cmtip]{xy}
\usepackage{kantlipsum}
\usepackage{tikz-cd}
\allowdisplaybreaks

\DeclareMathOperator{\SN}{SN}
\DeclareMathOperator{\id}{id}
\DeclareMathOperator{\pr}{pr}
\DeclareMathOperator{\Der}{Der}
\DeclareMathOperator{\Bott}{Bott}
\DeclareMathOperator{\Ha}{H}
\DeclareMathOperator{\Td}{Td}
\DeclareMathOperator{\tr}{tr}
\DeclareMathOperator{\str}{str}

\makeatletter
 \def\title@font{\normalsize\bfseries}
 \let\ltx@maketitle\@maketitle
 \def\@maketitle{\bgroup%
 \let\ltx@title\@title%
 \def\@\title{\resizebox{\textwidth}{!}{%
  \mbox{\title@font\ltx@title}%
 }}%
 \ltx@maketitle%
 \egroup}
\makeatother

\newcommand{\abs}[1]{\lvert#1\rvert}

\newcommand{\T}{\mathfrak{T}}
\newcommand{\A}{\mathcal{A}}

\newcommand{\R}{\mathbb{R}}
\newcommand{\Z}{\mathbb{Z}}

\newcommand{\g}{\mathfrak{g}}

\newcommand{\bd}{\begin{displaymath}}
\newcommand{\ed}{\end{displaymath}}
\newcommand{\be}{\begin{equation}}
\newcommand{\ee}{\end{equation}}

\numberwithin{equation}{section}
\theoremstyle{plain}
\newtheorem{corollary}[equation]{Corollary}

\newtheorem{Thm}{Theorem}
\newtheorem{lemma}[equation]{Lemma}
\newtheorem{proposition}[equation]{Proposition}

\theoremstyle{definition}

\theoremstyle{remark}
\newtheorem{remark}[equation]{Remark}

\begin{document}
\def\bp{\begin{proof}}
\def\ep{\end{proof}}
\def\C{\mathbb{C}}
\def\CE{\mathrm{CE}}
\def\D{\mathcal{D}}
\def\Dol{\mathrm{Dol}}
\def\E{\mathcal{E}}
\def\F{\mathcal{F}}
\def\H{\textbf{H}}
\def\k{\mathbb{K}}
\def\G{\mathcal{G}}
\def\M{\mathcal{M}}
\def\N{\mathcal{N}}
\def\m{\mathfrak{m}}
\def\t{\mathfrak{t}}
\def\O{\mathcal{O}}
\def\P{\mathcal{P}}
\def\U{\mathcal{U}}
\def\V{\mathscr{V}}
\def\L{\mathcal{L}}
\def\W{\mathcal{W}}
\def\X{\mathcal{X}}
\def\Y{\mathcal{Y}}
\def\ZZ{\mathcal{Z}}
\def\Im{\operatorname{Im}}
\def\Ber{\operatorname{Ber}}
\def\End{\operatorname{End}}
\def\Hom{\operatorname{Hom}}
\def\sgn{\operatorname{sgn}}
\def\rk{\operatorname{rank}}
\def\sh{\operatorname{sh}}

\title{Atiyah and Todd classes arising from integrable distributions}

\author{Zhuo Chen}
\address{Department of Mathematics, Tsinghua University}
\email{\href{mailto:~~~ zchen@math.tsinghua.edu.cn}{zchen@math.tsinghua.edu.cn}}

\author{Maosong Xiang}
\address{Beijing International Center for Mathematical Research, Peking University}
\email{\href{mailto:~~~ msxiang@pku.edu.cn}{msxiang@pku.edu.cn}}

\author{Ping Xu}
\thanks{Research partially supported by  NSF grants DMS-1406668 and DMS-1707545}
\address{Department of Mathematics, Pennsylvania State University}
\email{\href{mailto:~~~ ping@math.psu.edu}{ping@math.psu.edu}}

\begin{abstract}
  In this paper, we study the Atiyah class and Todd class of the DG manifold $(F[1],d_F)$ corresponding to an integrable distribution $F \subset T_\k M = TM \otimes_\R \k$, where $\k = \R$ or $\C$. We show that these two classes are canonically identical to those of the Lie pair $(T_\k M, F)$. As a consequence, the Atiyah class of a complex manifold $X$ is isomorphic to the Atiyah class of the corresponding DG manifold $(T^{0,1}_X[1],\bar{\partial})$. Moreover, if $X$ is a compact K\"{a}hler manifold, then the Todd class of $X$ is also isomorphic to the Todd class of the corresponding DG manifold $(T^{0,1}_X[1],\bar{\partial})$.
\end{abstract}

\maketitle

\tableofcontents

\section{Introduction}

Atiyah classes were introduced by Atiyah \cite{Atiyah} to characterize the obstruction to the existence of holomorphic connections on a holomorphic vector bundle over a complex manifold. In
\cite{Kap}, Kapranov showed that the Atiyah class $\alpha_X \in
 H_{}^1(X,T_X^\vee \otimes \End(T_X))$ of
 a   K\"{a}hler  manifold $X$ induces an $L_\infty$ algebra structure on $\Omega^{0,\bullet-1}_X(T_X)$ (see \cite{LaurentSX-CR, LaurentSX} for the complex
manifold case),
 which plays an important role in the formalism of Rozansky-Witten theory. In his seminar work \cite{Kon}, Kontsevich indicated a deep link between the Todd genus of complex manifolds and the Duflo element of Lie algebras.
See \cite{LSX} for the formality theorem for smooth DG manifolds,
which implies the Kontsevich-Duflo theorem for Lie algebras and
 Kontsevich's theorem for complex manifolds~\cite{Kon} under a
 unified framework.
 Inspired by the work of Kapranov and Kontsevich, Chen, Sti\'{e}non and Xu introduced the Atiyah class of Lie algebroid pairs in \cite{CSX}, which encodes both the Atiyah class of complex manifolds and the Molino class~\cite{Molino1} of foliations as special cases. Towards a different direction, Mehta, Sti\'{e}non and Xu~\cite{MSX} introduced Atiyah classes and Todd classes of DG manifolds. The Duflo element of a Lie algebra $\g$ is identified with the Todd class of the associated DG manifold $(\g[1],d_{\CE})$.


Let $F \subset T_\k M = TM \otimes_\R \k$ be a subbundle whose space of sections forms an integrable distribution, where $\k$ is either the field $\R$ or $\C$. When $\k = \R$, it corresponds to the tangent bundle of a regular foliation $\mathcal{F}$ in $M$. We will also say that $F$ is a regular foliation. 
Then the graded manifold $F[1]$ together with the Chevalley-Eilenberg differential $d_F$ on the Lie algebroid $F$ over $\k$ is a DG manifold. Thus we may consider the Atiyah class and the Todd class~\cite{MSX} of this DG manifold. Meanwhile, $(T_\k M,F)$ is a Lie algebroid pair over $M$. We may also consider the Atiyah class and the Todd class~\cite{CSX} of the Lie pair $(T_\k M,F)$.

The main purpose of this paper is to investigate the relation between these two types of Atiyah classes and Todd classes. First of all, we prove the following theorem in Section~\ref{Sec:ProofofThem1} regarding the cohomology groups of tensor fields of the DG manifold $(F[1],d_F)$.
\begin{Thm}\label{canonical isomorphism}
For any pair of non-negative integers $(n,m)$, there exists a canonical isomorphism
   \bd
     \Phi_m^{n}:~~~  \Ha^\bullet(\T_m^{n}(T(F[1])),L_{d_F}) \xrightarrow{\cong} \Ha_{\CE}^\bullet(F,\mathbb{T}_m^{n}(T_\k M/F))
   \ed
  from the cohomology of $(n,m)$-tensor fields $\T_m^{n}(T(F[1]))$ of the DG manifold $(F[1],d_F)$ to the Chevelley-Eilenberg cohomology of the $F$-module $\mathbb{T}_m^{n}(T_\k M/F)$ (see Equation~\eqref{notation} for notations).
  Moreover, the following diagram commutes:~~~
  \bd
  \xymatrix{
  \Ha^\bullet(\T_{m_1}^{n_1}(T(F[1])),L_{d_F}) \otimes_\k \Ha^\bullet(\T_{m_2}^{n_2}(T(F[1])),L_{d_F}) \ar[r]^-{\otimes} \ar[d]^{(\Phi_{m_1}^{n_1},\Phi_{m_2}^{n_2})} & \Ha^\bullet(\T_{m_1+m_2}^{n_1+n_2}(T(F[1])),L_{d_F}) \ar[d]^{\Phi_{m_1+m_2}^{n_1+n_2}} \\
  \Ha^\bullet_{\CE}(F,\mathbb{T}_{m_1}^{n_1}(T_\k M/F)) \otimes_\k \Ha_{\CE}^\bullet(F,\mathbb{T}_{m_2}^{n_2}(T_\k M/F)) \ar[r]^-{\otimes} & \Ha_{\CE}^\bullet(F,\mathbb{T}_{m_1+m_2}^{n_1+n_2}(T_\k M/F)).
  }
  \ed
\end{Thm}

Moreover, we have a homotopy contraction between $\Gamma(M,\wedge^\bullet F^\vee \otimes \wedge^\bullet (T_\k M/F))$ and the space of polyvector fields $\Gamma(F[1],\wedge^\bullet T(F[1]))$ of the DG manifold $(F[1],d_F)$.  Note that
$$
(\Gamma(F[1],\wedge^\bullet T(F[1])), L_{d_F}, \wedge, [-,-]_{SN})
$$
is a differential Gerstenhaber algebra, thus also a $(-1)$-shifted derived Poisson algebra~\cite{BCSX}, where $[-,-]_{SN}$ is the Schouten-Nijenhuis bracket. Op cit, a $(-1)$-shifted derived Poisson algebra structure is constructed on $\Gamma(M,\wedge^\bullet F^\vee \otimes \wedge^\bullet (T_\k M/F))$. We show in Proposition~\ref{application in Tpoly} that this structure in fact results from the homotopy transfer of the differential Gerstenhaber algebra structure in $\Gamma(F[1],\wedge^\bullet T(F[1]))$ via the aforementioned contraction data.

Then we prove the following main theorem in Section~\ref{Proofsection}. 
\begin{Thm}\label{Atiyah et todd classes}
 Under the canonical isomorphism as in Theorem~\ref{canonical isomorphism},
\begin{enumerate}
 \item The Atiyah class $\alpha_{F[1]}$ of the DG manifold $(F[1],d_F)$ corresponds to the Atiyah class $\alpha_{T_\k M/F}$ of the Lie pair $(T_\k M,F)$;
 \item The scalar Atiyah classes $c_k(F[1])$ ($k = 1,2,\cdots$) of the DG manifold $(F[1],d_{F})$ correspond  to the scalar Atiyah classes $c_k(T_\k M/F)  $ of the Lie pair $(T_\k M, F)$;
 \item The Todd class $\Td_{F[1]}$ of the DG manifold $(F[1],d_F)$ corresponds to the Todd class $\Td_{T_\k M/F} $ of the Lie pair $(T_\k M,F)$.
 \end{enumerate}
\end{Thm}

As a main application, consider a complex manifold $X$. Let $T_X^{1,0}$ and $T_X^{0,1}$ be its holomorphic and anti-holomorphic tangent bundle, respectively.
By considering the particular integrable
 distribution $F := T^{0,1}_X\subset T_\C X$, we get a Lie pair $(T_\C X, T_X^{0,1})$ and a DG manifold $(T_X^{0,1}[1], \bar{\partial})$, where $\bar{\partial}:\;\Omega^{0,\bullet}_X \rightarrow \Omega^{0,\bullet+1}_X$ is the Dolbeault differential operator. Applying Theorems~\ref{canonical isomorphism} and~\ref{Atiyah et todd classes}, we have the following
\begin{Thm}
Let $\Theta_X$ and $\Omega_X$ be the holomorphic tangent sheaf and holomorphic cotangent sheaf of $X$, respectively. Then
  \begin{enumerate}
    \item There exists a canonical isomorphism
            \bd
               \Phi_{m}^{n}:~~~  \Ha^\bullet(\mathfrak{T}_{m}^{n}(T(T_X^{0,1}[1])), L_{\bar{\partial}}) \xrightarrow{\cong} \Ha^\bullet(X,(\Omega_X)^{\otimes m} \otimes_{\O_X} (\Theta_X)^{\otimes n}),
            \ed
          where $\Ha^\bullet(X,(\Omega_X)^{\otimes m} \otimes_{\O_X} (\Theta_X)^{\otimes n})$ denotes the sheaf cohomology of $(\Omega_X)^{\otimes m} \otimes_{\O_X} (\Theta_X)^{\otimes n}$.
    \item The canonical isomorphism $\Phi_2^{1}$ sends the Atiyah class $\alpha_{T_X^{0,1}[1]}$ of the DG manifold $(T_X^{0,1}[1],\bar{\partial})$ to the Atiyah class $\alpha_X \in \Ha^1(X, \Omega_X^{\otimes 2} \otimes_{\O_X} \Theta_X)$ of the complex manifold $X$.
    \item The canonical isomorphism $\Phi_{\bullet}^{0}$ sends the Todd class $\Td_{T_X^{0,1}[1]}$ of the DG manifold $(T_X^{0,1}[1],\bar{\partial})$ to the Todd class $\Td_{T_X^{1,0}} \in \oplus_{k}\Ha^{k,k}(X)$ of the Lie pair $(T_\C X,T_X^{0,1})$. 
  \end{enumerate}
\end{Thm}

\textbf{Notations}:

In this note, the symbol $\k$ denotes either the field $\R$ or $\C$, and graded means $\Z$-graded. By a graded manifold $\M$ over $\k$, we mean a pair $(M,\O_\M)$, where $M$ is a smooth manifold, and $\O_\M$ is a sheaf of graded commutative $\k$-algebras over $M$, locally isomorphic to $C^\infty(U,\k) \otimes_\k \widehat{S(V^\vee)}$ for a graded vector space $V$ over $\k$, and $U \subset M$ open. We denote by $C^\infty(\M,\k)$ the space of $\k$-valued smooth functions on $\M$, i.e., the space of global sections of $\O_\M$. A graded vector bundle is a vector bundle in the category of graded manifolds(cf.\cite{mehta}).

For any graded vector bundle $\E$ over a graded manifold $\M$, denote by
\be\label{notation}
 \mathbb{T}_p^{q}(\E) = \E^{\otimes q} \otimes (\E^\vee)^{\otimes p}, \quad\quad\quad \T_p^q(\E) = \Gamma(\M;\mathbb{T}_p^q(\E)),
\ee
the tensor product of $q$-copies of $\E$ and $p$-copies of its dual $\E^\vee$, and the corresponding global section space, respectively. In particular, $\T_0^1(T\M) = \Gamma(T\M)$ is the space of vector fields of $\M$, which will also be denoted by $\mathfrak{X}(\M)$.

A DG manifold is a pair $(\M,Q)$, where $\M$ is a graded manifold, and $Q \in \mathfrak{X}(\M)$ is a homological vector field on $\M$, i.e., a degree $+1$ derivation of $C^\infty(\M)$ such that $[Q,Q] = 0$.
A DG vector bundle is a vector bundle in the category of DG manifolds(cf.\cite{mehta}).

Let $\E \rightarrow \M$ be a graded vector bundle and assume that $\M$ is a DG manifold. Then a DG-module structure on $\Gamma(\E)$ over $C^\infty(\M)$ induces a homological vector field on $\E$ such that $\E \rightarrow \M$ is a DG vector bundle. In fact, the category of DG vector bundles and the category of locally free DG-modules are equivalent(cf.~\cite{MSX}).

Let $(\M,Q)$ be a DG manifold. Then for any $p,q \geq 0$, the vector bundle
$$
\mathbb{T}_p^{q}(T\M) = (T\M)^{\otimes q} \otimes (T^\vee\M)^{\otimes p} ,
$$
together with the Lie derivative $L_{Q}$ along the homological vector field $Q$, is a DG vector bundle over $(\M,Q)$. And
$
(\wedge^q (T\M) \otimes \wedge^p (T^\vee\M), L_Q)
$
is a DG subbundle of $(\mathbb{T}_p^{q}(T\M),L_{Q})$.

\textbf{Acknowledgements:~~~ }
Xiang is grateful to his advisor Xiaobo Liu for his constant support and encouragement, to Pennsylvania State University for its hospitality and to China Scholarship Council for the financial support during his 20-months stay at Penn State. We would like to thank Ruggero Bandiera for telling us reference~\cite{Manetti} on the tensor product of homotopy contractions. We also wish to thank Damien Broka, Hsuan-Yi Liao and Mathieu Sti\'{e}non for helpful discussions.

\section{Cohomology of DG manifolds out of integrable distributions}\label{Sec:ProofofThem1}
In this section, we prove Theorem~\ref{canonical isomorphism} concerning the cohomology of $(n,m)$-tensor fields on the DG manifold $(F[1],d_F)$ via constructing a homotopy contraction explicitly.

\subsection{The contraction data}
Let $F \subset T_\k M$ be an integrable distribution. For simplicity, we will denote the quotient bundle $T_\k M/F$ by $B$ from now on. There is a short exact sequence of vector bundles over $M$:~~~
\be\label{SES II}
\xymatrix{
0 \ar[r] & F \ar[r]^-{i} & T_\k M \ar[r]^-{\pr_B} & B = T_\k M/F \ar[r] & 0.
}
\ee
There exists a flat $F$-connection $\nabla^{\Bott}$ on $B$, called Bott representation, which is defined by
\bd
 \nabla^{\Bott}_V Z:= \pr_{B}[V,j(Z)],\;\;\;\forall V \in \Gamma(F), Z \in \Gamma(B),
\ed
where $j: B \rightarrow T_\k M$ is any splitting of the short exact sequence~\eqref{SES II}.

Meanwhile, there associates to $F$ a DG manifold $(F[1], d_F)$, whose space of $\k$-valued smooth functions is
$$
  C^\infty(F[1],\k) = \Gamma(M;\wedge^\bullet F^\vee) :~= \Omega_F^\bullet.
$$
Here $d_F: \Omega_F^\bullet \rightarrow \Omega_{F}^{\bullet+1}$ is the Chevalley-Eilenberg differential of $F$. Let $\pi: F[1] \rightarrow M$ be the bundle projection map. For any vector bundle $E$ over $M$, we denote the space of smooth sections of the pull-back bundle $\pi^\ast(E)$ over $F[1]$ by
\bd
\Gamma(F[1],\pi^\ast(E)) \cong \Gamma(M;\wedge^\bullet F^\vee \otimes E) := \Omega_F^\bullet(E).
\ed
In particular, the covariant derivative $d_B$ of the Bott connection $\nabla^{\Bott}$ on $B$ gives rise to a cochain complex $(\Omega_F^\bullet(B),d_B)$, which is called the Chevalley-Eilenberg complex of $F$ with respect to $B$. And $d_B$ is also called the $F$-module structure on $B$.

The following lemma appeared without explicit proof in~\cite{Luca}. For completeness, we give a proof here.
\begin{lemma}\label{DR of VF}
  For any splitting $j:\;B \rightarrow T_\k M$ of the short exact sequence~\eqref{SES II}, there exists a contraction data
 \bd
  \begin{tikzcd}
   (\mathfrak{X}(F[1]),L_{d_F}) \arrow[loop left, distance=2em, start anchor={[yshift=-1ex]west}, end anchor= {[yshift=1ex]west}]{}{h_j} \arrow[r,yshift = 0.7ex, "\phi"] & (\Omega^\bullet_F(B),d_B) \arrow[l,yshift = -0.7ex, "\psi_j"],
  \end{tikzcd}
 \ed
  satisfying
  \begin{align}
  \label{dreq1} \phi \circ \psi_j &= \id:~~~  (\Omega^\bullet_F(B),d_B) \rightarrow (\Omega^\bullet_F(B),d_B), \\
  \label{dreq2} \id - \psi_j \circ \phi &= [L_{d_F},h_j]:~~~  (\mathfrak{X}(F[1]),L_{d_F}) \rightarrow (\mathfrak{X}(F[1]),L_{d_F}),
 \end{align}
 with side conditions
 \begin{align*}
   h_j \circ \psi_j &= 0, & \phi \circ h_j &= 0, & h_j^2 &= 0,
 \end{align*}
 where 
 $d_B$ is the Chevalley-Eilenberg differential of $F$ with respect to the $F$-module $B$.
\end{lemma}
\bp
Let us choose a splitting of the short exact sequence \eqref{SES II}, i.e., a pair of bundle maps $j:~~~  B \rightarrow T_\k M$ and $\pr_F:~~~  T_\k M \rightarrow F$ such that $\pr_B \circ j = \id_{B},\; \pr_F \circ i = \id_F$, and $\id_{T_\k M} = j \circ \pr_B + i \circ \pr_F$:~~~
\bd
 \xymatrix{
 0 \ar[r] & F \ar@<.7ex>[r]^-{i} & T_\k M \ar@<.7ex>[l]^-{\pr_F} \ar@<.7ex>[r]^-{\pr_B} & B \ar@<.7ex>[l]^-{j} \ar[r] & 0.
 }
\ed
By abuse of notations, we will use same notations to denote the corresponding splitting of the short exact sequence: 
\bd
 \xymatrix{
 0 \ar[r] & \Omega^\bullet_F(F) \ar@<.7ex>[r]^-{i} & \Omega^\bullet_F(T_\k M) \ar@<.7ex>[l]^-{\pr_F} \ar@<.7ex>[r]^-{\pr_B} & \Omega^\bullet_F(B) \ar@<.7ex>[l]^-{j} \ar[r] & 0.
 }
\ed
In the sequel, we consider $B$ as a subbundle of $T_\k M$ transversal to $F$ via the isomorphism $T_\k M \cong F \oplus B$ determined by $j$. The symbol $R$ denotes the space $C^\infty(M,\k)$ of $\k$-valued smooth functions on $M$. The space of vector fields on $F[1]$ is canonically isomorphic to the graded Lie algebra of $\k$-linear derivations of $\Omega^\bullet_F$:~~~
$$
\mathfrak{X}(F[1]) \cong \Der(\Omega^\bullet_F).
$$

First of all, we define the maps $\phi,\psi_j$ and $h_j$ as follows:~~~

\begin{itemize}
  \item Let $\phi$ be the $\Omega^\bullet_F$-linear map
  \be\label{phi}
   \phi:~~~ = \pr_B \circ \pi_\ast:~~~  \mathfrak{X}(F[1]) \rightarrow \Omega^\bullet_F(B),
  \ee
  where $\pi_\ast:~~~  \mathfrak{X}(F[1]) \rightarrow \Omega_F^\bullet(T_\k M)$ is the tangent map of the bundle projection $\pi:~~~  F[1] \rightarrow M$. In fact, $\pi_\ast(\X) \in \Omega_F^\bullet(T_\k M) \cong \Omega^\bullet_F \otimes_R \Der(R)$ is determined by
  \bd
   \pi_\ast(\X)(f) = \X(\pi^\ast f) \in \Omega_F^\bullet,\;\;\qquad\forall f \in R.
  \ed

  \item The map
  \begin{align}\label{psij}
  \psi_j:~~~  \Omega^\bullet_F(B) &\rightarrow \mathfrak{X}(F[1])
  \end{align}
  is the $\Omega^\bullet_F$-linear extension of the $R$-linear map
  \bd
   \psi_j:~~~  \Gamma(B) \rightarrow \Der(\Omega^\bullet_F) \cong \mathfrak{X}(F[1]) 
  \ed
  specified by the following relations:
  \begin{align}\label{defofpsij}
   \psi_j(Z)(\pi^\ast f) &= j(Z)(f), &  \psi_j(Z)(\xi) &= \pr_{F^\vee}(L_{j(Z)} \xi) \in \Omega_F^1,
  \end{align}
  for all $Z \in \Gamma(B), f \in R$ and $\xi \in \Omega^1_F$. More precisely, we have
  \be\label{defpsijv}
   \iota_V(\psi_j(Z)(\xi)) = j(Z)\langle V, \xi \rangle - \langle \pr_F[j(Z),V], \xi \rangle,
  \ee
  for all $V \in \Gamma(F)$. It is easy to see that the map $\psi_j$ is well-defined.

  \item Let $h_j:~\mathfrak{X}(F[1]) \rightarrow \mathfrak{X}(F[1])$ be the $\Omega^\bullet_F$-linear map of degree $(-1)$ determined by
   \be\label{h}
    h_j(\X) = (-1)^{\abs{\X}}\iota_{\pr_F(\pi_\ast(\X))} \in \Der(\Omega^\bullet_F) \cong \mathfrak{X}(F[1])
   \ee
   for all $\X \in \mathfrak{X}(F[1])$. Here $\pr_F(\pi_\ast(\X)) \in \Omega^\bullet_F(F)$, and $\iota$ is the contraction operator
  \bd
   \iota_{\alpha \otimes_R V}\omega := \alpha \wedge \iota_V\omega, \;\;\forall \alpha,\omega \in \Omega_F^\bullet, V \in \Gamma(F).
  \ed
  In particular, we have, for any $f \in R$,
  \begin{align}\label{equonh}
    h_j(\X)(\pi^\ast f) &= 0, & h_j(\X)(d_F(\pi^\ast f)) &= (-1)^{\abs{\X}}\pr_F(\pi_\ast(\X))(f).
  \end{align}
  \end{itemize}
The rest of proof is divided into several steps:~~~
\begin{itemize}
 \item \emph{The map $\phi$ is a cochain map.}


  It suffices to show that for any $\X \in \mathfrak{X}(F[1])$,
  \bd
   \phi([d_F,\X]) = d_B(\phi(\X)) \in \Omega^\bullet_F(B).
  \ed
  We prove, for any $\omega \in \Gamma(B^\vee)$, that the identity
  \be\label{keytochainmapofphi}
   \langle \omega, \phi([d_F,\X]) \rangle = \langle \omega, d_B(\phi(\X)) \rangle \in \Omega^\bullet_F
  \ee
  holds.
Note that $\pr_B^\vee(\omega) \in \Omega^1(M)$ is locally
of the form $\sum_k g_kdf_k$, where $f_k, g_k \in R$.
 It thus follows that
  \be\label{Eq1:phiiscochain}
  \sum g_k d_F(\pi^\ast f_k) = i^\vee\left(\sum g_kdf_k\right) = i^\vee(\pr_B^\vee(\omega)) = 0.
  \ee
  Now
\begin{align*}
  \text{LHS of}~\eqref{keytochainmapofphi} &= \langle \pr_B^\vee(\omega), \pi_\ast[d_F,\X] \rangle = \sum g_k \langle df_k, \pi_\ast[d_F,\X] \rangle = \sum g_k[d_F,\X](\pi^\ast f_k) \\
  &= \sum g_k d_F(\X(\pi^\ast f_k)) - \sum (-1)^{\abs{\X}}g_k\X(d_F\pi^\ast f_k) \quad\text{by Equation~\eqref{Eq1:phiiscochain}}\\
  &= d_F\left(\sum g_k\X(\pi^\ast f_k)\right) - \sum d_F(\pi^\ast g_k)\X(\pi^\ast f_k) + (-1)^{\abs{\X}}\sum \X(\pi^\ast g_k)d_F(\pi^\ast f_k)\\
  &= d_F\langle\omega, \phi(\X)\rangle - \sum d_F(\pi^\ast g_k)\X(\pi^\ast f_k) + \sum d_F(\pi^\ast f_k)\X(\pi^\ast g_k).
  \end{align*}
  Here in the last equality we have used the following identity
  \bd
   \langle\omega, \phi(\X)\rangle = \langle \pr_B^\vee(\omega),\pi_\ast(\X) \rangle = \sum g_k \langle df_k, \pi_\ast(\X) \rangle = \sum g_k \X(\pi^\ast f_k).
  \ed
  On the other hand, to compute the right-hand side of Equation~\eqref{keytochainmapofphi}, we first note that the
Lie algebroid cohomology differential $d_{B^\vee}$
corresponding to the dual $F$-module structure  on $B^\vee$
 can be computed by the following commutative diagram
  \bd
  \xymatrix{
   \Gamma(B^\vee) \ar[rr]^-{d_{B^\vee}} \ar[d]^{\pr_B^\vee} && \Omega_F^1(B^\vee) \subset \Omega_F^1(T_\k M)\\
   \Omega^1(M) \ar[rr]^-{d_{\text{DR}}} && \Omega^2(M). \ar[u]^-{R_F} 
  }
  \ed
  Here the right vertical map $R_F$ is defined by
  \bd
   R_F(\alpha \wedge \beta) = (i^\vee \otimes \id)(\alpha \otimes \beta - \beta \otimes \alpha) = i^\vee(\alpha)\otimes\beta - i^\vee(\beta)\otimes\alpha,
  \ed
  for all $\alpha, \beta \in \Omega^1(M)$.
  And the fact that $\Im(R_F \circ d_{\text{DR}} \circ \pr_B^\vee) = \Omega_F^1(B^\vee) \subset \Omega_F^1(T^\vee_\k M)$ follows from the assumption that $F \subset T_\k M$ is integrable. Thus,
  \begin{align*}
    \text{RHS of}~\eqref{keytochainmapofphi} &= d_F \langle \omega, \phi(\X) \rangle - \langle d_{B^\vee}(\omega), \phi(\X) \rangle \\
    &= d_F \langle \omega, \phi(\X) \rangle - \langle R_F\left(\sum dg_k\wedge df_k\right), \pr_B(\pi_\ast(\X)) \rangle \\
    &= d_F \langle \omega, \phi(\X) \rangle - \left\langle \pr_B^\vee\left(\sum d_F(\pi^\ast g_k) \otimes df_k - d_F(\pi^\ast f_k) \otimes dg_k\right), \pi_\ast(\X) \right\rangle \\
    &= d_F \langle \omega, \phi(\X) \rangle - \sum d_F(\pi^\ast g_k)\X(\pi^\ast f_k) + \sum d_F(\pi^\ast f_k)\X(\pi^\ast g_k).
  \end{align*}

  \item \emph{The map $\psi_j$ is a cochain map}, i.e., $\psi_j \circ d_B = L_{d_F} \circ \psi_j$. We have
    \begin{align*}
      \psi_j (d_B(\alpha \otimes_R Z)) &= d_F(\alpha) \otimes_R \psi_j(Z) + (-1)^{\abs{\alpha}}\alpha \otimes_R \psi_j(d_B(Z))
    \end{align*}
    and
    \begin{align*}
      L_{d_F}(\psi_j(\alpha \otimes_R Z)) &= d_F(\alpha) \otimes_R \psi_j(Z) + (-1)^{\abs{\alpha}}\alpha \otimes_R L_{d_F}(\psi_j(Z)),
   \end{align*}
   for any $\alpha \otimes_R Z \in \Omega^\bullet_F(B)$, where $\alpha \in \Omega^\bullet_F, Z \in \Gamma(B)$. Thus it suffices to prove that
   $$
    \psi_j(d_B(Z)) = L_{d_F}(\psi_j(Z)) = [d_F,\psi_j(Z)] \in \mathfrak{X}(F[1]).
   $$
  Let us examine this relation on the generators $\{\pi^\ast f,d_F(\pi^\ast f):~~~  f \in R\}$ of the $\k$-algebra $\Omega_F^\bullet$:

  On the one hand, for any $V \in \Gamma(F)$, by the defining Equation~\eqref{defofpsij} of $\psi_j$, we have
  \begin{align*}
  \iota_V(\psi_j(d_B(Z))(\pi^\ast f)) &= \psi_j(\nabla^{\Bott}_V(Z))(\pi^\ast f) = j(\pr_B([V,j(Z)]))(f) \\
   &= [V,j(Z)](f) - \pr_F([V,j(Z)])(f) \\
   &= V(j(Z)(f)) - (j(Z)(V(f)) - \pr_F([j(Z),V])(f))\;\quad\quad\text{by Equation~\eqref{defpsijv}}\\
   &= \iota_V([d_F,\psi_j(Z)](\pi^\ast f)).
  \end{align*}
  Therefore, we have
  \be\label{psijonf}
  \psi_j(d_B(Z))(\pi^\ast f) = [d_F,\psi_j(Z)](\pi^\ast f).
  \ee
  On the other hand, we have
  \begin{align*}
  [d_F,\psi_j(Z)](d_F(\pi^\ast f)) &= -d_F([d_F,\psi_j(Z)](\pi^\ast f))\qquad\qquad\qquad\qquad\qquad\text{by Equation\eqref{psijonf}}  \\
                        &= -d_F(\psi_j(d_B(Z))(\pi^\ast f))\\
                        &= -[d_F,\psi_j(d_B(Z))](\pi^\ast f) + \psi_j(d_B(Z))(d_F(\pi^\ast f)) \qquad\text{by Equation\eqref{psijonf}}  \\
                        &= -\psi_j(d^2_B(Z))(\pi^\ast f) + \psi_j(d_B(Z))(d_F(\pi^\ast f))\\
                        &= \psi_j(d_B(Z))(d_F(\pi^\ast f)).
  \end{align*}
  Hence, $\psi_j$ is a cochain map.

  \item \emph{Proof of Equation~\eqref{dreq1}}.
  We first claim that
  \be\label{dpiandpsij}
    \pi_\ast \circ \psi_j = j:~~~  \quad \Omega^\bullet_F(B) \rightarrow \Omega^\bullet_F(T_\k M),
  \ee
  Since both $\psi_j$ and $\pi_\ast$ are $\Omega^\bullet_F$-linear, it suffices to show that
  \bd
   \pi_\ast(\psi_j(Z)) = j(Z) \in \Gamma(T_\k M) \cong \Der(R),\;\;\forall Z \in \Gamma(B).
  \ed
  In fact, by the defining formula~\eqref{defofpsij} of $\psi_j$, we have
  \begin{align*}
    \pi_\ast(\psi_j(Z))(f) &= \psi_j(Z)(\pi^\ast f) = j(Z)(f),\;\;\forall f \in R.
  \end{align*}
  This proves Equation~\eqref{dpiandpsij}. Hence,
   \begin{align*}
   \phi \circ \psi_j &= \pr_B \circ \pi_\ast \circ \psi_j = \pr_B \circ j = \id.
   \end{align*}

\item \emph{Proof of Equation~\eqref{dreq2}}.
  We prove that
  \bd
   \X - \psi_j(\phi(\X)) = [L_{d_F},h_j](\X) = [d_F,h_j(\X)] + h_j([d_F,\X]),
  \ed
  for all $\X \in \mathfrak{X}(F[1])$. It suffices to show that both sides coincide when acting on elements $\{\pi^\ast f,d_F(\pi^\ast f): f \in R\}$. 
  That is, for any $f \in R$,
  \begin{align}\label{eq1}
      (\X - \psi_j(\Phi(\X)))(\pi^\ast f) &= [d_F,h_j(\X)](\pi^\ast f) + h_j([d_F,\X])(\pi^\ast f) = [d_F,h_j(\X)](\pi^\ast f),\\
  (\X - \psi_j(\Phi(\X)))(d_F(\pi^\ast f)) &= [d_F,h_j(\X)](d_F(\pi^\ast f)) + h_j([d_F,\X])(d_F(\pi^\ast f)). \label{eq2}
  \end{align}
  In fact, Equation~\eqref{eq1} can be verified directly:
  \begin{align*}
  [d_F,h_j(\X)](\pi^\ast f) &= (-1)^{\abs{\X}}h_j(\X)(d_F(\pi^\ast f)) \qquad\qquad\qquad\qquad\qquad\quad\text{by Equation~\eqref{equonh}} \\
              &=\pr_F(\pi_\ast(\X))(f) = (\pi_\ast(X) - \pr_B(\pi_\ast(\X)))(f)\;\quad\text{by Equation~\eqref{phi}} \\
              &=(\pi_\ast(\X) - (\phi(\X)))(f) \\
              &= \X(f) - \psi_j(\phi(X))(f).
  \end{align*}
   As for Equation~\eqref{eq2}, we compute
  \begin{align*}
   &\quad (\X - \psi_j(\phi(\X)))(d_F(\pi^\ast f)) \\
   &= (-1)^{\abs{\X}}(d_F((\X - \psi_j(\phi(\X)))(\pi^\ast f)) - [d_F, \X - \psi_j(\phi(\X))](\pi^\ast f))\qquad\quad\qquad\text{by Equation~\eqref{eq1}} \\
   &= (-1)^{\abs{\X}}(d_F([d_F,h_j(\X)](\pi^\ast f)) - [d_F,\X](\pi^\ast f) + [d_F,\psi_j(\phi(\X))](\pi^\ast f)) \\
   &= (-1)^{\abs{\X}}(d_F([d_F,h_j(\X)](\pi^\ast f)) - [d_F,\X](\pi^\ast f) + \psi_j(\phi([d_F,\X]))(\pi^\ast f)) \\ 
   &= [d_F,h_j(\X)](d_F(\pi^\ast f)) - (-1)^{\abs{\X}}(\pi_\ast([d_F,\X])(f) - j(\pr_B(\pi_\ast([d_F,\X])))(f)) \\
   &= [d_F,h_j(\X)](d_F(\pi^\ast f)) + (-1)^{\abs{\X}+1}\pr_F(\pi_\ast([d_F,\X]))(f) \qquad\qquad\qquad\qquad\;\text{by Equation~\eqref{equonh}}\\ 
   &= [d_F,h_j(\X)](d_F(\pi^\ast f)) + h_j([d_F,\X])(d_F(\pi^\ast f)).
  \end{align*}

  \item \emph{Side conditions}. Note that by Equation~\eqref{dpiandpsij},
   \bd
    (h_j \circ \psi_j)(Z) = \iota_{\pr_{F}(\pi_\ast(\psi_j(Z)))} = \iota_{\pr_F(j(Z))} = 0,\;\;\forall Z \in \Gamma(B).
   \ed
   Since both $h_j$ and $\psi_j$ are $\Omega_F^\bullet$-linear, it follows that $h_j \circ \psi_j = 0$. The verifications of $\phi \circ h_j = 0$ and $h_j^2 = 0$ are similar.
\end{itemize}
This completes the proof.
\ep

We generalize this result to the space of tensor fields on $F[1]$ as follows:~~~
\begin{proposition}\label{DR}
  For any splitting $j:~ B \rightarrow T_\k M$ of the short exact sequence~\eqref{SES II}, there exists a contraction data
  \bd
  \begin{tikzcd}
  (\T_m^{n}(T(F[1])),L_{d_F}) \arrow[loop left, distance=2.0em, start anchor={[yshift=-1ex]west}, end anchor={[yshift=1ex]west}]{}{H_{m}^n} \arrow[r,yshift = 0.7ex, "\Phi_m^{n}"]  & (\Omega^\bullet_F(\mathbb{T}_m^{n}(B)),d_B) \arrow[l,yshift = -0.7ex, "\Psi_m^{n}"],
  \end{tikzcd}
  \ed
 satisfying
 \begin{align*}
   \Phi_m^{n} \circ \Psi_m^{n} &= \id, & \Psi_m^{n} \circ \Phi_m^{n} + [L_{d_F},H_m^{n}] &= \id,
 \end{align*}
 with side conditions
 \begin{align*}
   H_m^{n} \circ \Psi_m^{n} &=0, & \Phi_m^{n} \circ H_m^{n} &= 0, & H_m^{n} \circ H_m^{n} &= 0.
 \end{align*}
 Here for simplicity we denote by $d_B$ the Chevalley-Eilenberg differential of $F$ with respect to the $F$-module $\mathbb{T}_m^{n}(B)$ for any pair of non-negative integers $(m,n)$.
\end{proposition}
\bp
  Note that all maps $\psi_j,\phi,h_j$ of the contraction data are $\Omega_F^\bullet$-linear. Their $\Omega_F^\bullet$-dual maps $\psi_j^\vee,\phi^\vee,h_j^\vee$ give rise to a contraction as well:~~~
  \bd
  \begin{tikzcd}
    (\Omega^1(F[1]),L_{d_F}) \arrow[loop left, distance=2em, start anchor={[yshift=-1ex]west}, end anchor={[yshift=1ex]west}]{}{h_j^\vee} \arrow[r,yshift = 0.7ex, "\psi_j^\vee"] & (\Omega^\bullet_F(B^\vee),d_{B}) \arrow[l,yshift = -0.7ex, "\phi^\vee"].
  \end{tikzcd}
 \ed
By considering the $\Omega_F^\bullet$-dual form of Equations~\eqref{dreq1} and~\eqref{dreq2}, it only suffices to show that
\bd
 [L_{d_F},h_j]^\vee = [L_{d_F},h_j^\vee]: \Omega^1(F[1]) \rightarrow \Omega^1(F[1]).
\ed
In fact, for all $\alpha \in \Omega^1(F[1]), \X \in \mathfrak{X}(F[1])$, we have
\begin{align*}
  \langle [L_{d_F},h_j^\vee]\alpha, & \X \rangle = \langle L_{d_F}(h_j^\vee(\alpha)), \X \rangle + \langle h_j^\vee(L_{d_F}(\alpha)), \X \rangle \\
  &= L_{d_F}\langle h_j^\vee(\alpha), \X \rangle - (-1)^{\abs{\alpha}-1}\langle h_j^\vee(\alpha), L_{d_F}\X \rangle + (-1)^{\abs{\alpha}+1}\langle L_{d_F}\alpha, h_j(\X) \rangle \\
  &= (-1)^{\abs{\alpha}}(L_{d_F}\langle \alpha, h_j(\X) \rangle - \langle L_{d_F}\alpha, h_j(\X) \rangle) + \langle \alpha, h_j(L_{d_F}\X) \rangle \\
  &= \langle \alpha, L_{d_F}(h_j(\X)) + h_j(L_{d_F}\X) \rangle = \langle \alpha, [L_{d_F},h_j]\X \rangle \\
  &= \langle [L_{d_F},h_j]^\vee(\alpha), \X \rangle.
\end{align*}

Define
  \begin{align}
  &\Psi_m^{n} :~~~ = (\phi^\vee)^{\otimes m} \otimes (\psi_j)^{\otimes n}:~~~  (\Omega^\bullet_F(\mathbb{T}_m^{n}(B)),d_B) \rightarrow (\T_m^{n}(T(F[1])),L_{d_F}), \label{psimn}\\
  &\Phi_m^{n} :~~~ = (\psi_j^\vee)^{\otimes m} \otimes (\phi)^{\otimes n}:~~~  (\T_m^{n}(T(F[1])),L_{d_F}) \rightarrow (\Omega^\bullet_F(\mathbb{T}_m^{n}(B)),d_B), \label{phimn}\\
  &H_m^{n}= \sum_{i=1}^m(\phi^\vee\circ\psi_j^\vee)^{\otimes(i-1)}\otimes h_j^\vee \otimes (\id)^{\otimes (n+m-i)} + \sum_{l=1}^n(\phi^\vee\circ\psi_j^\vee)^{\otimes m}\otimes (\psi_j\circ\phi)^{\otimes (l-1)} \otimes h_j \otimes \id^{\otimes (n-l)}. \label{Hmn}
  \end{align}
Then according to Lemma~\ref{DR of VF}, and the tensor trick in~\cite{Manetti}*{Example 2.6}, the maps $\Psi_m^{n}, \Phi_m^{n}, H_m^{n}$ constructed above satisfy the contraction properties.
\ep

Moreover, the quasi-isomorphisms $\Psi_m^{n}$ and $\Phi_m^{n}$ are, in fact, canonical up to homotopy:
\begin{proposition}\label{contraction and splitting}
 Assume that $\widehat{j}:~B \rightarrow T_\k M$ is another splitting of the short exact sequence~\eqref{SES II}. Let
 \bd
  \begin{tikzcd}
  (\T_m^{n}(T(F[1])),L_{d_F}) \arrow[loop left, distance=2.0em, start anchor={[yshift=-1ex]west}, end anchor={[yshift=1ex]west}]{}{\widehat{H}_m^n} \arrow[r,yshift = 0.7ex, "\widehat{\Phi}_m^n"]  & (\Omega^\bullet_F(\mathbb{T}_m^{n}(B)),d_B) \arrow[l,yshift = -0.7ex, "\widehat{\Psi}_m^n"]
  \end{tikzcd}
  \ed
  be the associated contraction. Then there exist two $\Omega^\bullet_F$-linear maps
 \begin{align*}
   \Theta_m^{n}:~~~  &\Omega^\bullet_F(\mathbb{T}_m^{n}(B)) \rightarrow \T_m^{n}(T(F[1])), &
   \Xi_m^{n}:~~~  &\T_m^{n}(T(F[1])) \rightarrow \Omega^\bullet_F(\mathbb{T}_m^{n}(B))
 \end{align*}
 such that
 \begin{align}
  \widehat{\Psi}^{n}_{m} - \Psi^{n}_m &= \Theta_m^{n} \circ d_B + L_{d_F} \circ \Theta_m^{n}:~~~  \Omega^\bullet_F(\mathbb{T}_m^{n}(B)) \rightarrow \T_m^{n}(T(F[1])), \label{ds1}\\
  \widehat{\Phi}_m^{n} - \Phi_m^{n} &= \Xi_m^{n} \circ L_{d_F} + d_B \circ \Xi_m^{n}:~~~  \T_m^{n}(T(F[1])) \rightarrow \Omega^\bullet_F(\mathbb{T}_m^{n}(B)). \label{ds2}
  \end{align}
\end{proposition}

According to the defining formulas~\eqref{psimn} and~\eqref{phimn} of $\Psi^{n}_m$ and $\Phi^{n}_m$ in the proof of Proposition~\ref{DR}, the key step is to study how the quasi-isomorphism $\psi_j$ depends on the splitting $j$:

\begin{lemma}\label{different splitting}
Under the same assumptions as in Proposition \ref{contraction and splitting}, there exists a $\Omega^\bullet_F$-linear
homotopy between $\widehat{\Psi}^1_0 = \psi_{\widehat{j}}$ and $\psi_j$
$$
\Theta:\;\Omega^\bullet_F(B)  \rightarrow \mathfrak{X}(F[1]),
$$
i.e.,
\begin{align}\label{homotopy retraction}
  \psi_{\widehat{j}} - \psi_j &= \Theta \circ d_B + L_{d_F} \circ \Theta:~~~  \Omega^\bullet_F(B) \rightarrow \mathfrak{X}(F[1]).
\end{align}
\end{lemma}
\bp
Let $\theta = \widehat{j} - j \in \Gamma(M;\Hom(B,F))$. Define a map
\bd
\Theta:~~~  \Omega^\bullet_F(B) \rightarrow \Der(\Omega_F^\bullet) \cong \mathfrak{X}(F[1])
\ed
by
\begin{align}\label{Theta}
  \Theta(\alpha \otimes_R Z) &= (-1)^{\abs{\alpha}}\alpha \otimes \iota_{\theta(Z)},
\end{align}
where $\alpha \in \Omega_F^\bullet, Z \in \Gamma(B)$.

We show that this $\Theta$ satisfies Equation~\eqref{homotopy retraction}. In fact, both sides of~\eqref{homotopy retraction} are $\Omega^\bullet_F$-linear, since
\begin{align*}
L_{d_F}(\Theta(\alpha \otimes_R Z))  &= (-1)^{\abs{\alpha}}[d_F,\alpha \otimes_R \Theta(Z)] = (-1)^{\abs{\alpha}}d_F(\alpha) \otimes_R \Theta(Z) + \alpha \otimes_R L_{d_F}(\Theta(Z)), \\
 \Theta(d_B(\alpha \otimes_R Z)) &= \Theta(d_F(\alpha) \otimes_R Z + (-1)^{\abs{\alpha}}\alpha \otimes_R d_B(Z)) \\
                                &= (-1)^{\abs{\alpha}+1}d_F(\alpha) \otimes_R \Theta(Z) + \alpha \otimes_R \Theta(d_B(Z)),
\end{align*}
for all $\alpha \otimes_R Z \in \Omega^\bullet_F(B)$, where $\alpha \in \Omega_F^\bullet, Z \in \Gamma(B)$.

It thus suffices to show that
\be\label{eq3}
\psi_{\widehat{j}}(Z) - \psi_j(Z) = \Theta(d_B(Z)) + L_{d_F}(\Theta(Z)) \in \mathfrak{X}(F[1]).
\ee
Let us denote by $\pr_F,\widehat{\pr}_F:~T_\k M \rightarrow F$ the projections corresponding to the two splittings $j$ and $\widehat{j}$, respectively.
Then
\begin{align*}
  \widehat{\pr}_F &= \pr_F - \theta \circ \pr_B:~~~  \Omega^\bullet_F(T_\k M) \rightarrow \Omega^\bullet_F(F).
\end{align*}
For any $f \in R$, 
\begin{align*}
   (\Theta(d_B(Z)) + L_{d_F}(\Theta(Z)))(\pi^\ast f) &= [d_F, \Theta(Z)](\pi^\ast f) = [d_F,\iota_{\theta(Z)}](\pi^\ast f) = \theta(Z)(\pi^\ast f) \\
   &= (\psi_{\widehat{j}}(Z) - \psi_j(Z))(\pi^\ast f); 
\end{align*}
For any $\xi \in \Omega^1_F, V \in \Gamma(F)$, we have
\begin{align*}
 &\quad \iota_V(\psi_{\widehat{j}}(Z)(\xi) - \psi_j(Z)(\xi)) \qquad\qquad\qquad\text{by Equation~\eqref{defpsijv}}\\
 &= \theta(Z)\langle V,\xi \rangle - \langle \widehat{\pr}_F([\widehat{j}(Z),V]) - \pr_F([j(Z),V]),\xi\rangle \\
 &= \theta(Z)\langle V,\xi \rangle - \langle\pr_F([\widehat{j}(Z),V]) - \theta(\pr_B[\widehat{j}(Z),V]) - \pr_F([j(Z),V]), \xi \rangle \\
 &= \theta(Z)\langle V,\xi \rangle - \langle [\theta(Z),V], \xi \rangle  - \langle \theta(\nabla^{\Bott}_VZ), \xi \rangle \\
 &= \langle V, L_{\theta(Z)}(\xi) \rangle - \langle \theta(\iota_Vd_B(Z)), \xi \rangle = \iota_{V} (L_{\theta(Z)}(\xi)) + \iota_V(\Theta(d_B(Z))(\xi))\\
 &= \iota_V(L_{d_F}(\Theta(Z))(\xi) + \Theta(d_B(Z))(\xi)).  
\end{align*}
Since $\Omega_F^\bullet$ is generated by $R$ and $\Omega_F^1$, it follows that Equation~\eqref{eq3} holds.
\ep

\bp[Proof of Proposition~\ref{contraction and splitting}]
Define
  \begin{align*}
    \Theta_m^{n} &= \sum_{l=1}^n(\phi^\vee)^{\otimes m}\otimes(\psi_{\widehat{j}})^{\otimes (l-1)}\otimes\Theta\otimes (\psi_j)^{\otimes(n-l)}:~~~  \Omega^\bullet_F(\mathbb{T}_m^{n}(B)) \rightarrow \T_m^{n}(T(F[1])), \\
    \Xi_m^{n} &= \sum_{k=1}^m(\psi_{\widehat{j}}^\vee)^{\otimes (k-1)}\otimes\Theta^\vee \otimes(\psi_j^\vee)^{\otimes(m-k)}\otimes(\phi)^{\otimes n}:~~~  \T_m^{n}(T(F[1])) \rightarrow \Omega^\bullet_F(\mathbb{T}_m^{n}(B)).
  \end{align*}
  Here $\Theta^\vee$ denotes the $\Omega^\bullet_F$-dual of $\Theta$ defined in Equation~\eqref{Theta}.

Then by Lemma~\ref{different splitting}, Equations~\eqref{ds1} and~\eqref{ds2} follow from straightforward computations.
\ep

\subsection{Proof of Theorem~\ref{canonical isomorphism}}

Now we are ready to prove Theorem~\ref{canonical isomorphism}:~~~
\bp[Proof of Theorem~\ref{canonical isomorphism}]
Let $j:\;B \rightarrow T_\k M$ be a splitting of the short exact sequence~\eqref{SES II}. Applying Proposition~\ref{DR}, we obtain an isomorphism
 \bd
  \Phi^{n}_m:~~~  \Ha^\bullet(\T_m^{n}(T(F[1])),L_{d_F}) \xrightarrow{\cong} \Ha_{\CE}^\bullet(F,\mathbb{T}_m^{n}(B)).
 \ed
 According to Proposition~\ref{contraction and splitting}, $\Phi_m^{n}$ does not depend on the splitting $j$, and thus is canonical.
 Finally, the commutative diagram in Theorem~\ref{canonical isomorphism} follows from the commutative diagram on the cochain level
 \bd
 \xymatrix{
 \T_{m_1}^{n_1}(T(F[1])) \otimes \T_{m_2}^{n_2}(T(F[1])) \ar[r]^-{\otimes} \ar[d]^{\Phi_{m_1}^{n_1}\otimes\Phi_{m_2}^{n_2}} & \T_{m_1+m_2}^{n_1+n_2}(T(F[1])) \ar[d]^{\Phi_{m_1+m_2}^{n_1+n_2}} \\
  \Omega^\bullet_F(\mathbb{T}_{m_1}^{n_1}(B)) \otimes \Omega^\bullet_F(\mathbb{T}_{m_2}^{n_2}(B)) \ar[r]^-{\otimes} & \Omega^\bullet_F(\mathbb{T}_{m_1+m_2}^{n_1+n_2}(B)).
 }
 \ed
 This completes the proof.
\ep

In the sequel, the contraction data $(\Phi_m^{n},\Psi_m^{n},H_m^{n})$ in Proposition~\ref{DR} will simply be denoted by $(\Phi,\Psi,H)$:~~~
\bd
  \begin{tikzcd}
  (\T_m^{n}(T(F[1])),L_{d_F}) \arrow[loop left, distance=2em, start anchor={[yshift=-1ex]west}, end anchor={[yshift=1ex]west}]{}{H}  \arrow[r,yshift = 0.7ex, "\Phi"] & (\Omega^\bullet_F(\mathbb{T}_m^{n}(B)),d_B) \arrow[l,yshift = -0.7ex, "\Psi"].
  \end{tikzcd}
\ed
It is clear that $(\Phi,\Psi,H)$ induces a contraction data on the subcomplexes as well
 \be\label{contractionforOmegaTpoly}
  \begin{tikzcd}
  (\Gamma(F[1], \wedge^nT(F[1]) \otimes \wedge^mT^\vee(F[1])), L_{d_F}) \arrow[loop left, distance=2em, start anchor={[yshift=-1ex]west}, end anchor={[yshift=1ex]west}]{}{H} \arrow[r,yshift = 0.7ex, "\Phi"] & (\Omega^\bullet_F(\wedge^n B \otimes \wedge^{m}B^\vee),d_B) \arrow[l,yshift = -0.7ex, "\Psi"].
  \end{tikzcd}
 \ee
The following result is an immediate consequence, which implies that cohomologies of $\Gamma(F[1], \wedge^nT(F[1]) \otimes \wedge^mT^\vee(F[1]))$ are bounded from above:~~~
\begin{corollary}\label{cohomology of forms}
 The cohomology
  \bd
   \Ha^r(\Gamma(F[1], \wedge^nT(F[1]) \otimes \wedge^mT^\vee(F[1])),L_{d_F})
  \ed
 vanishes under any of the following conditions:~~~
 \begin{enumerate}
   \item  $m > \rk(B)$; \item $n > \rk(B)$;
   \item $r > \rk(F)$.
 \end{enumerate}
\end{corollary}
\subsection{Application}
Consider the case $m=0$. The contraction data reduces to
 \be\label{contractionforTpoly}
  \begin{tikzcd}
  ( \Gamma(F[1],\wedge^{n}T(F[1])),L_{d_F}) \arrow[loop left, distance=2em, start anchor={[yshift=-1ex]west}, end anchor={[yshift=1ex]west}]{}{H} \arrow[r,yshift = 0.7ex, "\Phi"] & (\Omega^\bullet_F(\wedge^n B),d_B) \arrow[l,yshift = -0.7ex, "\Psi"],
  \end{tikzcd}
 \ee
where
\begin{align*}
  \Psi = \wedge^n\psi_j:\; &\Omega^\bullet_F(\wedge^n B) \rightarrow \Gamma(F[1],\wedge^{n}T(F[1])), \\
  \Phi = \wedge^n\phi:\; &\Gamma(F[1],\wedge^{n}T(F[1])) \rightarrow \Omega^\bullet_F(\wedge^n B), \\
   H = \sum_{l=1}^n \wedge^{(l-1)} (\psi_j \circ \phi) \wedge h_j \wedge (\id)^{\wedge(n-l)}: \; &\Gamma(F[1],\wedge^{n}T(F[1])) \rightarrow \Gamma(F[1],\wedge^{n}T(F[1])).
\end{align*}
Note that both $\Omega^\bullet_F(\wedge^\bullet B)$ and $\Gamma(F[1],\wedge^{\bullet}T(F[1]))$, together with the obvious wedge product, are commutative graded algebras. It is clear that both $\Psi$ and $\Phi$ are compatible with the wedge product, thus are morphisms of commutative algebras. However, $H$ in general is not a derivation. In fact, we have
\bd
 H(\X_1\wedge\X_2) = H(\X_1)\wedge \X_2 + (\Psi \circ \Phi)(\X_1) \wedge H(\X_2),\;\;\forall \X_i \in \Gamma(\wedge^\bullet T(F[1])).
\ed

Note that $(\Gamma(F[1],\wedge^{\bullet}T(F[1])),L_{d_F})$ together with the wedge product and the Schouten-Nijenhuis bracket $[-,-]_{\SN}$, forms a differential Gerstenhaber algebra. By the standard homotopy transfer technique, we recover the following
\begin{proposition}\label{application in Tpoly}Given a splitting $j:~~~  B \rightarrow T_\k M$ of the short exact sequence~\eqref{SES II}, then
\begin{enumerate}
  \item there exists a $(-1)$-shifted derived Poisson algebra structure $\{\lambda_k\}_{k \geq 1}$ on $\Omega_F^\bullet(\wedge^\bullet B)$ such that $\lambda_1 = d_B$. 
  \item this $(-1)$-shifted derived Poisson algebra structure coincides with the one out of the Lie pair $(T_\k M, F)$ as in~\cite{BCSX}*{Proposition 4.3}.
\end{enumerate}
\end{proposition}

The following technical lemma is needed:
\begin{lemma}\label{techniquelemma in Tpoly}
  For any $\X_1,\X_2 \in \mathfrak{X}(F[1])$ and $Z \in \Gamma(B)$, we have
  \begin{align}
    \phi([h_j(\X_1),h_j(\X_2)]) &= 0, \label{eq4}\\
    h_j([h_j(\X_1),\psi_j(Z)]) &= 0.   \label{eq5}
  \end{align}
\end{lemma}
\bp
Note that the space of vertical vector fields
$$
\mathfrak{X}(F[1])^{\text{vertical}} = \{\iota_{\mathcal{V}} \in \Der(\Omega_F^\bullet):~~~  \mathcal{V} \in \Omega^\bullet_F(F)\} \cong \ker(\pi_\ast)
$$
is a Lie subalgebra of $\mathfrak{X}(F[1])$. Since $\Im(h_j) \subset \mathfrak{X}(F[1])^{\text{vertical}}$, it follows that
\bd
 \phi([h_j(\X_1),h_j(\X_2)]) = \pr_B(\pi_\ast([h_j(\X_1),h_j(\X_2)])) = 0.
\ed

Meanwhile, Since
\bd
  \pi_\ast([h_j(\X_1),\psi_j(Z)])(f) = [h_j(\X_1),\psi_j(Z)](\pi^\ast f) = h_j(\X_1)(j(Z)(f)) = 0,
\ed
for all $f \in R$, it follows that $\pi_\ast([h_j(\X_1),\psi_j(Z)]) = 0$. According to Equation~\eqref{h}, the latter implies Equation~\eqref{eq5}, as desired.
\ep

\bp[Proof of Proposition~\ref{application in Tpoly}]
We first use the contraction data~\eqref{contractionforTpoly} to prove $(1)$:~~~  Note that both $\Omega_F^\bullet(\wedge^\bullet B)$ and $\Gamma(F[1],\wedge^\bullet T(F[1]))$ carry natural graded commutative algebra structures, i.e., the wedge products. Moreover,
$$
(\Gamma(F[1],\wedge^\bullet T(F[1])),L_{d_F},\wedge,[-,-]_{\SN})
$$
is a differential Gerstenhaber algebra, and therefore is also a $(-1)$-shifted derived Poisson algebra. Since both $\Psi$ and $\Phi$ are morphisms of the underlying commutative algebras, by the homotopy transfer for derived Poisson algebras~\cite{BCSX}*{Theorem 2.16}, there induces a $(-1)$-shifted derived Poisson algebra structure $\{\lambda_k\}_{k \geq 1}$ on $\Omega^\bullet_F(\wedge^\bullet B)$, whose first bracket $\lambda_1$ is $d_B$ and higher brackets are defined iteratively op cit.

To prove $(2)$, let us work out higher brackets $\{\lambda_k\}_{k \geq 2}$ explicitly:~~~
\begin{itemize}
  \item[$k=2$.] The binary bracket $\lambda_2$ is, by definition
    \bd
      \lambda_2(\ZZ_1,\ZZ_2) = \Phi([\Psi(\ZZ_1),\Psi(\ZZ_2)]_{\SN}),\;\forall \ZZ_1,\ZZ_2 \in \Omega^\bullet_F(\wedge^\bullet B).
    \ed
    More precisely, we consider the situation on generating elements:~~~  
  \begin{align*}
    \lambda_2(Z_1, Z_2) &= \phi([\psi_j(Z_1),\psi_j(Z_2)]_{\SN}) = \pr_B([j(Z_1),j(Z_2)]),\;\forall Z_i \in \Gamma(B), i =1,2, \\
    \lambda_2(Z_1,\xi) &= \phi([\psi_j(Z_1),\psi_j(\xi)]_{\SN}) = \psi_j(Z_1)(\xi) = \pr_{F^\vee}(L_{j(Z_1)}\xi),\;\forall Z_1 \in \Gamma(B), \xi \in \Omega_F^1.
  \end{align*}
  Similarly, it follows that
  \begin{align*}
    \lambda_2(Z_1,f) &= j(Z_1)(f), & \lambda_2(\omega_1,\omega_2) &= 0,
  \end{align*}
  for all $f \in R = C^\infty(M,\k)$ and $\omega_i \in \Omega^\bullet_F$.
  \item[$k=3$.] The trinary bracket $\lambda_3$ is given by
  \bd
   \lambda_3(\ZZ_1,\ZZ_2,\ZZ_3) = \sum_{\sigma \in \sh(2,1)}\pm \Phi([H([\Psi(\ZZ_{\sigma(1)}),\Psi(\ZZ_{\sigma(2)})]_{\SN}), \Psi(\ZZ_{\sigma(3)})]_{\SN}),
  \ed
  for all $\ZZ_i \in \Omega^\bullet_F(\wedge^\bullet B), 1 \leq i \leq 3$, where $\sh(2,1)$ is the set of $(2,1)$-shuffles, and $\pm$ is some proper Koszul sign. In particular, we have, on the generating elements,
  \bd
   \lambda_3(Z_1,Z_2,\xi) = \phi([h_j([\psi_j(Z_1),\psi_j(Z_2)]_{\SN}),\psi_j(\xi)]_{\SN}) = \iota_{\pr_F[j(Z_1),j(Z_2)]}(\xi),
  \ed
  for all $Z_1,Z_2 \in \Gamma(B)$ and $\xi \in \Gamma(F^\vee)$, where $\pr_F:~~~  \Gamma(T_\k M) \rightarrow \Gamma(F)$ is the projection. Furthermore, since
  \be\label{eq6}
  h_j([\X,\psi_j(\xi)]_{\SN}) = h_j(\X(\xi)) = 0,\;\;\;\forall \X \in \mathfrak{X}(F[1]),\xi \in \Gamma(F^\vee) \subset \Omega_F^\bullet,
  \ee
  it follows that $\lambda_3$ vanishes if restricted to $\Gamma(B) \times \Gamma(B) \times \Gamma(B), \Gamma(B) \times \Gamma(F^\vee) \times \Gamma(F^\vee)$ and $\Gamma(F^\vee) \times \Gamma(F^\vee) \times \Gamma(F^\vee)$.

  \item[$k=4$.] The 4th bracket $\lambda_4$ acting on the generating elements is of form
  \begin{align*}
   \lambda_4(Z_1,Z_2,Z_3,Z_4) &= \sum_{\sigma \in \sh(2,2)}\pm \phi([h_j([\psi_j(Z_{\sigma(1)}),\psi_j(Z_{\sigma(2)})]_{\SN}), h_j([\psi_j(Z_{\sigma(1)}), \psi_j(Z_{\sigma(2)})])]_{\SN}) \\
   &+ \sum_{\tau \in \sh(2,1,1)}\pm \phi([h_j([h_j([\psi_j(Z_{\tau(1)}),\psi_j(Z_{\tau(2)})]_{\SN}), \psi_j(Z_{\tau(3)})]_{\SN}),\psi_j(Z_{\tau(4)})]_{\SN}),
  \end{align*}
  for all $Z_i \in \Gamma(B), 1 \leq i \leq 4$. It follows from Lemma~\ref{techniquelemma in Tpoly} that $\lambda_4(Z_1,Z_2,Z_3,Z_4) = 0$. Moreover, it follows from Equations~\eqref{eq5} and~\eqref{eq6} that $\lambda_4$ also vanishes with one or more inputs from $\Gamma(F^\vee)$.
  \item[$k \geq 5$.] The $L_\infty$-brackets $\{\lambda_k\}_{k \geq 5}$ vanish for similar reasons.
\end{itemize}
It thus follows that the generating relations of the $(-1)$-shifted derived Poisson algebra structure on $\Omega^\bullet_F(\wedge^\bullet B)$, are exactly those of \cite{BCSX}*{Proposition 4.3}. This concludes the proof.
\ep

\begin{remark}
The $(-1)$-shifted derived Poisson algebra structure on $\Omega_F^\bullet(\wedge^\bullet B)$ induces a strong homotopy Lie-Rinehart algebra structure on the subspace $\Omega^\bullet_F(B)$. This result is due to Vitagliano~\cite{Luca}.
\end{remark}

According to Theorem 1.1 in~\cite{BCSX}, this $(-1)$-shifted derived Poisson algebra is in fact canonical up to isomorphism. On the other hand, $(\Ha^\bullet(\Gamma(F[1],\wedge^{\bullet}T(F[1])),L_{d_F}),\wedge, [-,-]_{\SN})$ is canonically a Gerstenhaber algebra. As an immediate consequence of Proposition~\ref{application in Tpoly} and Theorem~\ref{canonical isomorphism}, we have
\begin{corollary}
  There is a canonical isomorphism of Gerstenhaber algebras
  \bd
  \Psi:~~~  (\Ha^\bullet_{\CE}(F,\wedge^\bullet B), \wedge, \lambda_2) \xrightarrow{\cong} (\Ha^\bullet(\Gamma(F[1],\wedge^{\bullet}T(F[1])),L_{d_F}),\wedge,[-,-]_{\SN}).
  \ed
\end{corollary}

\section{Atiyah and Todd classes}\label{Proofsection}
\subsection{Atiyah and Todd classes of DG manifolds}
Let $(\M,Q)$ be a DG manifold. The Atiyah class $\alpha_\M \in \Ha^1(\T_2^{1}(T\M),L_Q)$ of $(\M,Q)$ is, by definition, the Atiyah class of the tangent DG vector bundle $T\M \rightarrow \M$ with respect to the DG Lie algebroid $T\M$ \cite{MSX}:~~~Choose a $T\M$-connection $\nabla$ on the tangent bundle $T\M$. Then there associates a degree $+1$ $L_Q$-cocycle $\alpha_\M^\nabla \in \T_2^{1}(T\M)$ defined by
\be\label{Atiyah cocycle}
 \alpha^\nabla_\M(X,Y) = L_Q(\nabla_X Y) - \nabla_{L_Q(X)}Y - (-1)^{\abs{X}}\nabla_X L_Q(Y),\;\forall X,Y \in \Gamma(T\M),
\ee
which will be called the Atiyah cocycle of $\M$ with respect to $\nabla$. The cohomology class
\bd
 \alpha_\M = [\alpha_\M^\nabla] \in \Ha^1(\T_2^{1}(T\M), L_Q),
\ed
which does not depend on the choice of $\nabla$, is called the Atiyah class of the DG manifold $(\M,Q)$.

Note that
$$
\alpha_\M \in \Ha^1(\T_2^{1}(T\M),L_Q) \cong \Ha^1(\Omega^1(\M) \otimes_{C^\infty(\M,\k)} \Gamma(\End(T\M)),L_Q).
$$
For any positive integer $k$, one can form $\alpha_\M^k$, the image of $\alpha_\M^{\otimes k}$ under the natural map
\bd
 \otimes_\k^k\Ha^1(\Omega^1(\M) \otimes_{C^\infty(\M,\k)} \Gamma(\End(T\M)),L_Q) \rightarrow \Ha^{k}(\Omega^k(\M) \otimes_{C^\infty(\M,\k)} \Gamma(\End(T\M)),L_Q) 
\ed
induced by the wedge product in $\Omega^\bullet(\M)$ and the composition in $\End(T\M)$.

The scalar Atiyah classes~\cite{MSX} of the DG manifold $(\M,Q)$ are defined by
\be\label{scalar atiyah classes}
 c_k(\M):~~~ = \frac{1}{k!}\left(\frac{i}{2\pi}\right)^k \str(\alpha_\M^k) \in \Ha^{k}(\Omega^k(\M),L_Q),\;k = 1,2,\cdots,
\ee
where $\str:~~~  \Gamma(\End(T\M)) \rightarrow C^\infty(\M,\k)$ is the supertrace map.

The Todd class~\cite{MSX} of the DG manifold $(\M,Q)$ is defined by
\be\label{Todd class}
  \Td_\M:~~~ = \Ber\left( \frac{\alpha_\M}{1 - e^{-\alpha_\M}} \right) \in \prod_{k \geq 0}\Ha^{k}(\Omega^k(\M),L_Q),
\ee
where $\Ber:~~~  \Gamma(\End(T\M)) \rightarrow C^\infty(\M,\k)$ is the Berezinian (or the superdeterminant) map~\cite{Manin}.

In this paper, we focus on the DG manifold $(F[1], d_F)$ corresponding to an integrable distribution $F \subset T_\k M$. The Atiyah, scalar Atiyah and Todd classes of $(F[1],d_F)$, are denoted respectively by:~~~
\begin{align*}
& \alpha_{F[1]} ~\in \Ha^1(\T_2^{1}(T(F[1])),L_{d_F}),
\\
& c_k(F[1]) ~\in \Ha^k(\Omega^k(F[1]),L_{d_F}), \quad k=1,2,\cdots
\\
& \Td_{F[1]} ~\in \prod_{k\geq 0}\Ha^k(\Omega^k(F[1]),L_{d_F}).
\end{align*}

\subsection{Atiyah and Todd classes of the Lie pair $(T_\k M,F)$}
Given an integrable distribution $F \subset T_\k M$, $(T_\k M,F)$ is a $\k$-Lie algebroid pair (or Lie pair for short). We briefly recall from~\cite{CSX} the Atiyah and Todd classes of the Lie pair $(T_\k M,F)$.

Recall that there is a short exact sequence of vector bundles over $M$:~~~
\bd
\xymatrix{
0 \ar[r] & F \ar[r] & T_\k M \ar[r] & T_\k M/F \ar[r] & 0.
}
\ed
Let us choose a splitting $j: T_\k M/F \rightarrow T_\k M$ of the above short exact sequence and a $T_\k M$-connection $\nabla$ on $T_\k M/F$ extending the Bott $F$-connection. The associated Atiyah cocycle
$$
R^{\nabla} \in \Omega^1_F(\mathbb{T}_2^{1}(T_\k M/F)) = \Gamma(M;F^\vee \otimes F^\perp \otimes \End(T_\k M/F)),
$$
where $F^\perp = (T_\k M/F)^\vee$, is defined by
\begin{align*}
 R^{\nabla}(V,Z)W &:= \nabla_{V}\nabla_{j(Z)}(W) - \nabla_{j(Z)}\nabla_{V} (W) - \nabla_{[V,j(Z)]}(W),
\end{align*}
for all $V \in \Gamma(F)$ and $Z,W \in \Gamma(T_\k M/F)$. The cohomology class
\bd
\alpha_{T_\k M/F} = [R^{\nabla}] \in \Ha_{\CE}^1(F,\mathbb{T}_2^{1}(T_\k M/F)) \cong \Ha_{\CE}^1(F,F^\perp \otimes \End(T_\k M/F))
\ed
does not depend on the choice of $j$ and $\nabla$, and is called the Atiyah class of the Lie pair $(T_\k M,F)$.
\begin{remark}
On the one hand, when $\k = \R$, this particular class is also known as the Molino class~\cite{Molino1} of the foliation $F$. On the other hand, the Atiyah class of the Lie pair $(T_\C X,T^{0,1}_X)$ corresponding to a complex manifold $X$ is identical to Atiyah class $\alpha_X$~\cite{Atiyah} of the complex manifold $X$.
\end{remark}
The scalar Atiyah classes of the Lie pair $(T_\k M,F)$ are defined by
\be\label{scalar atiyah class of pair}
 c_k(T_\k M/F):~~~ = \frac{1}{k!}\left(\frac{i}{2\pi}\right)^k \operatorname{tr}(\alpha^k_{T_\k M/F}) \in \Ha_{\CE}^k(F,\wedge^k F^\perp),\quad k = 1,2,\cdots.
\ee
Here $\alpha^k_{T_\k M/F}$ is the image of $\alpha^{\otimes k}_{T_\k M/F}$ under the natural map
\bd
 \otimes_\k^k(\Ha_{\CE}^1(F, F^\perp \otimes \End(T_\k M/F))) \rightarrow \Ha_{\CE}^k(F,\wedge^k F^\perp \otimes \End(T_\k M/F))
\ed
induced by the composition in $\End(T_\k M/F)$ and the wedge product in $\wedge^\bullet F^\perp$.

The Todd class of the Lie pair $(T_\k M,F)$ is the cohomology class
\be\label{Todd class of pair}
 \Td_{T_\k M/F} = \det\left(\frac{\alpha_{T_\k M/F}}{1 - e^{-\alpha_{T_\k M/F}}}\right) \in \bigoplus_{k\geq 0}\Ha_{\CE}^k(F,\wedge^kF^\perp).
\ee

\begin{remark}
  Let $X$ be a compact K\"{a}hler manifold. Then the natural inclusion
  $$
  \oplus_{k}\Ha^k(X,\Omega_X^k) \hookrightarrow \oplus_{k}\Ha^{2k}(X,\C)
  $$
  maps the scalar Atiyah classes $c_k(T^{1,0}_X)$ and the Todd class $\Td_{T_{X}^{1,0}}$ of the Lie pair $(T_\C X,T_X^{0,1})$ to the $k$-th
Chern characters $\operatorname{ch}_k(X)$ and the Todd class
 $\Td_X$ of $X$, respectively.
\end{remark}

\subsection{Splittings and connections}
From now on, we fix a splitting $j:~~~  B \rightarrow T_\k M $ of the short exact sequence~\eqref{SES II} so that we can identify $T_\k M$ with $F \oplus B$.

Let us choose a $T(F[1])$-connection $\nabla$ on $T(F[1])$. It induces an operator
\bd
 \nabla^B:~~~  \Omega^\bullet_F(B) \times \Omega^\bullet_F(B) \rightarrow \Omega^\bullet_F(B)
\ed
as follows:
\be\label{Eqforconnection}
 \nabla^B_\ZZ \W = \Phi(\nabla_{\Psi(\ZZ)}\Psi(\W)),\;\;\;\forall \ZZ,\W \in \Omega^\bullet_F(B).
\ee
Here $\Phi$ and $\Psi$ are part of the contraction data in the previous section. The following lemma can be easily verified:~~~
\begin{lemma}\label{lembeijing}
  \begin{itemize}
  \item The restriction of $\nabla^B$ defines a bilinear map
  \bd
   \nabla^B:~~~  \Gamma(B) \times \Gamma(B) \rightarrow \Gamma(B),
  \ed
  which together with the Bott-$F$-connection on $B$ determines a unique $T_\k M$-connection on $B$. This particular $T_\k M$-connection on $B$ will also be denoted by $\nabla^B$ by abuse of notation.
  \item For all $Z \in \Gamma(B)$ and $V \in \Gamma(F)$, we have
  \begin{align}\label{Leibniz rule}
    \iota_{\pr_F([Z,V])} = [\nabla^B_Z,\iota_V] = \nabla^B_Z \circ \iota_V - \iota_V \circ \nabla^B_Z:~~~  \quad \Omega_F^\bullet(B) \rightarrow \Omega_F^{\bullet-1}(B),
  \end{align}
  where $\iota_V:~~~  \Omega^\bullet_F(B) \rightarrow \Omega_F^{\bullet-1}(B)$ is the contraction operator, and $\pr_F:~~~  T_\k M \rightarrow F$ is the projection determined by the splitting $j$.
\end{itemize}
\end{lemma}

\subsection{Proof of Theorem~\ref{Atiyah et todd classes}}
\subsubsection{Atiyah classes}
By Theorem~\ref{canonical isomorphism}, the first statement of Theorem~\ref{Atiyah et todd classes} follows directly from the following Lemma:~~~
\begin{lemma}\label{lemonF}
  Under the same hypothesis as in Lemma~\ref{lembeijing}, the quasi-isomorphism
  $$
  \Phi = \Phi_2^{1}:~~~  \T_2^{1}(T(F[1])) \rightarrow \Omega_F^\bullet(\mathbb{T}_2^{1}(B))
  $$
  maps the Atiyah cocycle $\alpha_{F[1]}^{\nabla}$ of the DG manifold $(F[1],d_F)$ with respect to the connection $\nabla$, to the Atiyah cocycle $R^{\nabla^{B}}$ of the Lie pair $(T_\k M,F)$ with respect to the $T_\k M$-connection $\nabla^{B}$.
\end{lemma}
\bp
Recall that the Atiyah cocycle $\alpha_{F[1]}^{\nabla} \in \T_2^{1}(T(F[1]))$ is of degree $1$. Thus $\Phi(\alpha_{F[1]}^{\nabla}) \in \Omega^1_F(\mathbb{T}_2^{1}(B))$. It suffices to show that
 \bd
  \Phi(\alpha_{F[1]}^\nabla)(Z,W) = R^{\nabla^{B}}(-,Z)W \in \Omega^1_F(B),
 \ed
for all $Z,W \in \Gamma(B)$. In fact, by Equation~\eqref{phimn},
\begin{align*}
&\quad\Phi\left(\alpha_{F[1]}^\nabla\right)(Z,W) = \Phi\left(\alpha_{F[1]}^\nabla(\Psi(Z),\Psi(W))\right) \qquad\qquad\qquad\qquad\qquad\text{by Equation~\eqref{Atiyah cocycle}}\\
&= \Phi\left(L_{d_F}(\nabla_{\Psi(Z)}\Psi(W)) - \nabla_{L_{d_F}(\Psi(Z))}\Psi(W) - \nabla_{\Psi(Z)}L_{d_F}(\Psi(W))\right) \\
&= d_B\left(\Phi(\nabla_{\Psi(Z)}\Psi(W))\right) - \Phi\left(\nabla_{\Psi(d_B(Z))}\Psi(W)\right) - \Phi\left(\nabla_{\Psi(Z)}(\Psi(d_B(W)))\right) \;\;\text{by Equation~\eqref{Eqforconnection}}\\
&= d_B\left(\nabla^B_{Z}W\right) - \nabla^B_{d_B(Z)}W - \nabla^B_{Z}d_B(W).
\end{align*}
Hence, for any $V \in \Gamma(F)$,
\begin{align*}
  &\quad\iota_V\left(\Phi\left(\alpha_{F[1]}^\nabla\right)(Z,W)\right) = \iota_V\left(d_B\left(\nabla^B_{Z}W\right) - \nabla^B_{d_B(Z)}W - \nabla^B_{Z}d_B(W)\right) \\
  &=\nabla^{B}_V\nabla^{B}_{Z}W - \nabla^{B}_{\nabla^{B}_VZ}W - \iota_V(\nabla^B_Zd_B(W)) \qquad\qquad\qquad\quad\text{by Equation~\eqref{Leibniz rule}} \\
  &=\nabla^{B}_V\nabla^{B}_{Z}W - \nabla^{B}_{\nabla^{B}_VZ}W - \nabla^{B}_{Z}\nabla^{B}_VW + \nabla^{B}_{\pr_F[Z,V]}W \\
  &= \nabla^{B}_V\nabla^{B}_{Z}W - \nabla^{B}_{Z}\nabla^{B}_VW - \nabla^{B}_{[V,Z]}W = R^{\nabla^{B}}(V,Z)W.
\end{align*}
Here we have used the definition of Bott-connection on $B$, which implies that
\bd
 -\nabla^B_VZ + \pr_F[Z,V] = -\pr_B[V,Z] - \pr_F[V,Z] = -[V,Z].
\ed
This completes the proof.
\ep

\subsubsection{Scalar Atiyah classes}
In this section, we prove the second statement of Theorem~\ref{Atiyah et todd classes} about scalar Atiyah classes. Consider the supertrace map
\begin{align*}
  \str:~~~  \Gamma(\End(T(F[1]))) &\rightarrow \Omega^\bullet_F
\end{align*}
and the trace map
\begin{align*}
  \tr:~~~  \Omega^\bullet_F(\End(B)) &\rightarrow \Omega^\bullet_F.
\end{align*}
They are both cochain maps, and $\Omega^\bullet_F$-linear. In an obvious manner, both maps can be extended to
\begin{align*}
  \str:~~~  \Omega^m(F[1]) \otimes_{\Omega^\bullet_F} \Gamma(\End(T(F[1]))) &\rightarrow \Omega^m(F[1]), & \tr:~~~  \Omega^\bullet_F(\wedge^m B^\vee \otimes \End(B)) &\rightarrow \Omega^\bullet_F(\wedge^m B^\vee),
\end{align*}
which are also cochain maps.

Recall that by Proposition~\ref{DR}, we have a contraction data
  \be\label{contraction}
  \begin{tikzcd}
   (\Omega^m(F[1]) \otimes_{\Omega^\bullet_F} \Gamma(\End(T(F[1]))),L_{d_F}) \arrow[loop left, distance=2em, start anchor={[yshift=-1ex]west}, end anchor={[yshift=1ex]west}]{}{H} \arrow[r,yshift = 0.7ex, "\Phi"] & (\Omega^\bullet_F(\wedge^m  B^\vee  \otimes \End(B)),d_B) \arrow[l,yshift = -0.7ex, "\Psi"].
  \end{tikzcd}
  \ee

\begin{lemma}\label{STRetTr}
  The cochain maps $\Phi \circ \str$ and $\tr \circ \Phi$ are homotopic with $\Phi \circ \str \circ H$ being the homotopy map, i.e.,
  \bd
    \Phi \circ \str - \tr \circ \Phi = d_{B} \circ (\Phi \circ \str \circ H) + (\Phi \circ \str \circ H) \circ L_{d_F}:~~~  \Omega^m(F[1]) \otimes_{\Omega^\bullet_F} \Gamma(\End(T(F[1]))) \rightarrow \Omega^\bullet_F(\wedge^m B^\vee).
  \ed
\end{lemma}
\bp
First of all, we prove the following
 \begin{align}\label{stronB}
   \Psi \circ \tr &= \str \circ \Psi:~~~  \Omega^\bullet_F(\wedge^m B^\vee \otimes \End(B)) \rightarrow \Omega^m(F[1]).
 \end{align}
 When $m=0$, Equation \eqref{stronB} reduces to
 \begin{align*}
   \tr &= \str \circ \Psi_1^1:~~~  \Omega^\bullet_F(\End(B)) \rightarrow \Omega^\bullet_F,
 \end{align*}
which is clearly true,
since all maps involved are $\Omega_F^\bullet$-linear and elements in $\Gamma(B)$ are purely even.
For any $\omega \otimes_{\Omega^\bullet_F} \lambda \in \Omega^\bullet_F(\wedge^m B^\vee \otimes \End(B))$, where $\omega \in \Omega_F^\bullet(\wedge^m B^\vee)$ and $\lambda \in \Omega_F^\bullet(\End(B))$,
we have
\begin{align*}
 \str(\Psi(\omega \otimes_{\Omega_F^\bullet} \lambda)) &= \str((\wedge^m\phi^\vee)(\omega) \otimes_{\Omega_F^\bullet} \Psi_1^1(\lambda)) = (\wedge^m\phi^\vee)(\omega) \otimes_{\Omega_F^\bullet} \str(\Psi_1^1(\lambda)) \\
 &= \wedge^m\phi^\vee(\omega) \otimes_{\Omega_F^\bullet} \tr(\lambda) = (\wedge^m\phi^\vee)(\tr(\omega \otimes_{\Omega^\bullet_F} \lambda)) = \Psi(\tr(\omega \otimes_{\Omega^\bullet_F} \lambda)).
\end{align*}
This proves~\eqref{stronB}.

For any $\Lambda \in \Omega^m(F[1]) \otimes_{\Omega^\bullet_F} \Gamma(\End(T(F[1])))$, by the contraction data~\eqref{contraction}, we have
\bd
 \Lambda = \Psi(\Phi(\Lambda)) + [L_{d_F},H](\Lambda).
\ed
Thus it follows that
\begin{align*}
   \Phi(\str(\Lambda)) &= \Phi(\str(\Psi(\Phi(\Lambda)))) + \Phi(\str([L_{d_F},H](\Lambda))) \qquad\quad \text{by Equation~\eqref{stronB}}\\
                       &= (\Phi \circ \Psi)(\tr(\Phi(\Lambda))) + (\Phi \circ \str)(L_{d_F}(H(\Lambda))) + (\Phi \circ \str)(H(L_{d_F}(\Lambda)))\\
                       &= \tr(\Phi(\Lambda)) + d_B(\Phi(\str(H(\Lambda)))) + \Phi(\str(H(L_{d_F}(\Lambda)))).
\end{align*}
Here we used the fact that both $\Phi$ and $\str$ are cochain maps.
\ep

As an immediate corollary, we have
\begin{corollary}\label{strandtr}
  The following diagram
  \bd
  \xymatrix{
    \Ha^\bullet(\Omega^m(F[1])\otimes_{\Omega_F^\bullet} \Gamma(\End(T(F[1]))),L_{d_F}) \ar[d]^-{\str} \ar[r]^-{\Phi} & \Ha_{\CE}^\bullet(F,\wedge^{m}B^\vee \otimes \End(B)) \ar[d]^-{\tr} \\
    \Ha^\bullet(\Omega^m(F[1]),L_{d_F})  \ar[r]^-{\Phi} & \Ha_{\CE}^\bullet(F,\wedge^m B^\vee)
  }
  \ed
  commutes, i.e.,
  \bd
   \Phi \circ \str = \tr \circ \Phi :~~~  \Ha^\bullet(\Omega^m(F[1])\otimes_{\Omega^\bullet_F} \Gamma(\End(T(F[1]))),L_{d_F}) \rightarrow \Ha_{\CE}^\bullet(F,\wedge^m B^\vee).
  \ed
\end{corollary}

\bp[Proof of Theorem~\ref{Atiyah et todd classes} $(2)$]
 Recall the defining equations~\eqref{scalar atiyah classes} and~\eqref{scalar atiyah class of pair} of the two types of scalar Atiyah classes. Then
 \begin{align*}
   \Phi(c_k(F[1])) &= \frac{1}{k!}\left(\frac{i}{2\pi}\right)^k \Phi(\str(\alpha_{F[1]}^k))\qquad\quad \text{by Corollary~\ref{strandtr}} \\
                   &= \frac{1}{k!}\left(\frac{i}{2\pi}\right)^k \tr(\Phi(\alpha_{F[1]}^k))\;\quad \text{since $\Phi$ is an algebra morphism} \\
                   &= \frac{1}{k!}\left(\frac{i}{2\pi}\right)^k \tr((\Phi(\alpha_{F[1]}))^k)\qquad\quad \text{by Theorem~\ref{Atiyah et todd classes} $(1)$} \\
                   &= \frac{1}{k!}\left(\frac{i}{2\pi}\right)^k \tr(\alpha_{T_\k M/F}^k) = c_k(T_\k M/F).
 \end{align*}
\ep

\subsubsection{Todd classes}
Before we prove the last statement of Theorem~\ref{Atiyah et todd classes}, we recall a general method calculating the determinant and Berezinian of formal power series of cohomology classes. Let
\bd
 P(x) = 1 + \sum_{i\geq 1} b_i x^i, 
\ed
be a formal power series, where $b_i \in \k$.
According to~\cite{Hirz}*{Lemma 1.2.2}, one can associate to $P(x)$ a sequence of polynomials $K_r(p_1,\cdots,p_r), r = 1,2,\cdots,$ of homogeneous weight $r$, called $m$-sequence. Here by saying that $K_r$ is of homogeneous weight $r$, we mean that $K_r$ is of the form
$$
K_r(p_1,\cdots,p_r) = \sum\limits_{n_1 + 2n_2 + \cdots + rn_r = r}c_{n_1\cdots n_r}(b_i)p_1^{n_1}\cdots p_r^{n_r}.
$$
It is well-known that
\begin{align*}
  \Ber(P(\A)) &= 1 + \sum_{r \geq 1} K_r(\str(\A),\cdots,\str(\A^r)),  \\
  \det(P(A)) &= 1 + \sum_{r \geq 1} K_r(\tr(A),\cdots,\tr(A^r))
\end{align*}
for any supermatrix $\A$ in a graded commutative $\k$-algebra, and any matrix $A$ in an ordinary commutative $\k$-algebra.

Let
$$
s := \min\{\rk(F),\rk(B)\}.
$$
It follows from Corollary~\ref{cohomology of forms} that for any
 $\omega \in \Ha^1(\T_2^{1}(T(F[1])),L_{d_F})$, we have
 $\omega^r = 0$ for all $r > s$ . Since $K_r$ is of homogeneous weight $r$, one has
 $$
 \Ber(P(\omega)) = 1 + \sum_{r= 1}^s K_r(\str(\omega),\cdots,\str(\omega^r)).
$$
Similarly, for any $\lambda \in \Ha^1_{\CE}(F,B^\vee\otimes\End(B))$, one has
$$
 \det(P(\lambda)) = 1 + \sum_{r = 1}^sK_r(\tr(\lambda),\cdots, \tr(\lambda^r)).
$$

\begin{lemma}\label{connectingQclass}
 With notations as above, we have, for any $\omega \in \Ha^1(\T_2^{1}(T(F[1])),L_{d_F})$
  \bd
    \Phi(\Ber(P(\omega))) = \det(P(\Phi(\omega))) \in \bigoplus_{k \geq 0}\Ha^k_{\CE}(F,\wedge^k B^\vee).
  \ed
\end{lemma}
\bp
 By Corollary~\ref{strandtr} and the fact that $\Phi$ is an algebra morphism, we have
\begin{align*}
  \Phi(\str(\omega^i)) &= \tr(\Phi(\omega^i)) = \tr((\Phi(\omega))^i),~~\;~\forall 1 \leq i \leq s.
\end{align*}
Therefore,
\begin{align*}
  \Phi(\Ber(P(\omega))) &= 1 + \sum_{r = 1}^sK_r(\Phi(\str(\omega)),\cdots,\Phi(\str(\omega^r))) \\
                        &= 1 + \sum_{r = 1}^sK_r(\tr{\Phi(\omega)},\cdots,\tr{(\Phi(\omega))^r}) = \det(P(\Phi(\omega))),
\end{align*}
as desired.
\ep

Now we prove the last statement of Theorem~\ref{Atiyah et todd classes}:~~~
\bp[Proof of Theorem~\ref{Atiyah et todd classes} $(3)$]
 We take the particular power series
$$
P(x) = \frac{x}{1 - e^{-x}}.
$$
By the defining formulas~\eqref{Todd class} and~\eqref{Todd class of pair} of Todd classes, we have
$$
\Td_{F[1]}=\Ber P(\alpha_{F[1]}),\quad \Td_{T_\k M/F}=\det P(\alpha_{T_\k M/F}).
$$

By Theorem~\ref{Atiyah et todd classes} $(1)$, we have
\bd
  \alpha_{T_\k M/F}=\Phi(\alpha_{F[1]}).
\ed

Hence,
\begin{align*}
\Phi(\Td_{F[1]}) &=\Phi(\Ber P(\alpha_{F[1]})) \qquad\quad\text{by Lemma~\ref{connectingQclass}}\\
                 &= \det P (\Phi(\alpha_{F[1]}))\qquad\quad\text{by Theorem~\ref{Atiyah et todd classes} $(1)$} \\
                 &=\det P(\alpha_{T_\k M/F})=\Td_{T_\k M/F}.\\
\end{align*}
\ep

\end{document}